\documentclass[10pt,psamsfonts]{amsart}
\usepackage{amsthm,amsmath,amsfonts,amssymb,amsthm,graphicx,epstopdf,color,pifont,subfigure,flafter,bm,amscd}
\usepackage{graphicx,overpic}
\usepackage[margin=3.6cm]{geometry}
\usepackage[colorlinks=true, linkcolor=blue, citecolor=red]{hyperref}

\def\R{\mathbb R}
\def\sss{\mathbb S}

\newtheorem {theorem} {Theorem}
\newtheorem {proposition} [theorem]{Proposition}

\usepackage{float}
\usepackage{tikz, wrapfig}
\tikzset{node distance=3cm, auto}
\allowdisplaybreaks

\begin{document}

\title[Global centers]
{Global centers of a class of\\ cubic polynomial differential systems}

\author[Jaume Llibre and Gabriel Rondón]{Jaume Llibre$^{1}$ and Gabriel Rondón$^{2}$}

\address{$^{1}$Departament de Matemàtiques, Edifici Cc, Universitat Autònoma de Barcelona, 08193 Bellaterra, Barcelona, Catalonia, Spain}

\address{$^{2}$S\~{a}o Paulo State University (Unesp), Institute of Biosciences, Humanities and
Exact Sciences. Rua C. Colombo, 2265, CEP 15054--000. S. J. Rio Preto, S\~ao Paulo, Brazil.}
 
\email{jaume.llibre@uab.cat}
\email{gabriel.rondon@unesp.br}

\thanks{ .}


\date{}

\begin{abstract}
A difficult classical problem in the qualitative theory of differential systems in the plane $\R^2$ is the center-focus problem, i.e. to distinguish between a focus and a center. Another difficult problem is to distinguish inside a family of centers the ones which are global. A global center is a center $p$ such that $\mathbb{R}^2\setminus\{p\}$ is filled with periodic orbits.

In this paper we classify the global centers of the family of real polynomial differential systems of degree $3$ that in complex notation write
$$
i\dot{w}=w-A_3\overline{w}^2-A_4w^3-A_5w^2\overline{w}-A_6w\overline{w}^2,
$$
where $w=x+iy$ and $A_k\in\mathbb{C}$ for $k=3,4,5,6$.
\end{abstract}

\maketitle

\section{Introduction}

The general cubic differential equations in complex notation having either a weak focus or a center at the origin of coordinates are
\begin{equation}\label{main}
i\dot{w}=w-A_1w^2-A_2w\overline{w}-A_3\overline{w}^2-A_4w^3-A_5w^2\overline{w}-A_6w\overline{w}^2-A_7\overline{w}^3,
\end{equation}
where $w=x+iy$ and $A_k\in\mathbb{C},$ $k=1,\cdots,7$. These differential systems have been considered in several articles, thus in the papers \cite{MR1660202, MR1426207, 5f1ebd1d-940d-311b-a97b-a2b7787a895c, MR1383940} the authors provide for some subclasses of these differential systems necessary and sufficient conditions in order that the  equilibrium point at the origin of coordinates be a center. In this paper we restrict  our attention to the subclass $A_1=A_2=A_7=0$. Then the complex differential equation \eqref{main} can be written as the following real cubic polynomial differential system 
\begin{equation}\label{main_eq}
\begin{array}{rl}
\dot{x}=&y+2a_1xy-a_2(x^2-y^2)-(b_2+c_2+d_2)x^3-(3b_1+c_1-d_1)x^2y\\
&+(3b_2-c_2-d_2)xy^2,\\
\dot{y}=&-x+a_1(x^2-y^2)+2a_2xy+(b_1+c_1+d_1)x^3-(3b_2+c_2-d_2)x^2y\\
&+(-3b_1+c_1+d_1)xy^2+(b_2-c_2+d_2)y^3,
\end{array}
\end{equation}
in the plane $\mathbb{R}^2$, where $A_3=a_1+ia_2,$ $A_4=b_1+ib_2,$ $A_5=c_1+ic_2$ and $A_6=d_1+id_2.$

The main goal of this paper is to classify the global centers of system \eqref{main_eq}. We recall that a \textit{center} of a differential equation in $\mathbb{R}^2$ is an equilibrium point $p$ having a neighbourhood $\mathcal{U}$ such that $\mathcal{U}\setminus\{p\}$ is filled with periodic orbits. In particular, a \textit{global center} is a center $p$ such that $\mathbb{R}^2\setminus\{p\}$ is filled with periodic orbits.

The classical problem of distinguishing between a focus and a center is a difficult one and it is one of the challenges of the theory of nonlinear differential systems in the plane $\R^2$. The rigorous notion of center appeared in the literature with the works of Poincaré \cite{18} and Dulac \cite{6}. 

The classification of all centers of the cubic polynomial differential systems is a problem unsolved now. It has been solved only for some subclasses of cubic polynomial differential systems, thus the authors of \cite{5f1ebd1d-940d-311b-a97b-a2b7787a895c} provided necessary and sufficient conditions to determine when the origin of the differential system \eqref{main_eq} is a center. Their result is the following.

\begin{theorem}\label{lemma_main}
The cubic polynomial differential systems \eqref{main_eq} have a center at the origin if, and only if, one of the following four sets of conditions is satisfied:
\begin{itemize}
\item[(i)] $c_2=0,$ $b_2d_1+d_2b_1=0,$ $3b_1-d_1=0;$
\item[(ii)] $c_2=0,$ $b_2d_1+d_2b_1=0,$  $F=a_2^2d_2^3-3a_2^2d_2d_1^2+6a_2a_1d_2^2d_1-2a_2a_1d_1^3-a_1^2d_2^3+3a_1^2d_2d_1^2=0;$
\item[(iii)] $c_1=c_2=0,$ $b_2-d_2=b_1+d_1=0;$
\item[(iv)] $c_2=0=d_1=d_2=0,$ $G=-a_2^2b_2^3+3a_2^2b_2b_1^2+6a_2a_1b_2^2b_1-2a_2a_1b_1^3+a_1^2b_2^3-3a_1^2b_2b_1^2.$
\end{itemize}
\end{theorem}

There are several works studying the global centers of different classes of polynomial differential systems, for example in \cite{8} and \cite{LLIBRE2021125281} the authors proved that polynomial differential systems of even degree cannot have global centers because such systems always have orbits coming from and going to infinity. However to classify all polynomial differential systems of odd degree having a global center is a very difficult problem. This last problem was proposed by Conti in \cite{5}, and in fact up to now only few partial results exist for certain families of polynomial differential systems of odd degree, for more details see, for instance, \cite{9,11,doi:10.1080/14689367.2023.2228737}.

Roughly speaking the Poincar\'e compactification, first, consists in identifying the plane $\R^2$ with the interior of unit closed disc $\mathbb D^2$ centered at the origin of coordinates, and the boundary of this disc, the circle $\mathbb S^1$, with the infinity of $\R^2$. After the polynomial differential system defined in $\R^2$ is extended analytically to the whole closed disc $\mathbb D^2$. In this way we can study the dynamics of the polynomial differential systems in a neighbourhood of the infinity. All the details on the Poincar\'e compactification can be found in \cite[Chapter 5]{MR2256001}.

\begin{figure}[h]
\begin{overpic}[scale=0.55]{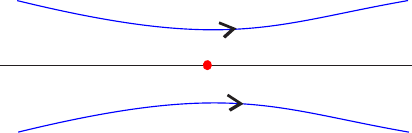}
\put(102,15){$\infty$}
\end{overpic}
\caption{\footnotesize{An equilibrium point on the infinity line whose local phase portrait is formed by two hyperbolic sectors whose two separatrices are contained on the infinity line.}}\label{fig_sechyp}
\end{figure}

We emphasize that in order that the origin of a cubic polynomial differential system \eqref{main_eq} be a global center, that system either does not have equilibrium points at infinity in the Poincar\'e compactification, or the local phase portraits of all its infinite equilibrium points are formed by two hyperbolic sectors having their two separatrices at infinity, see Figure \ref{fig_sechyp}. This implies that the Jacobian matrix at any infinite equilibrium point of system \eqref{main_eq} must be identically zero, otherwise the infinite equilibrium point would be either hyperbolic, or semi-hyperbolic, or nilpotent, and it is known that the local phase portraits of such kind of equilibrium points are not formed by two hyperbolic sectors having their two separatrices contained at infinity, for details see Theorems 2.15, 2.19 and 3.5 of \cite{MR2256001}. 	

Using Theorem \ref{lemma_main} we prove the main result of this article, stated in the next theorem.

\begin{theorem}\label{main_teo}
The polynomial differential systems \eqref{main_eq} have a global center at the origin of coordinates if, and only if, one of the following sets of conditions holds:
\begin{itemize}
\item[(a)] $a_1=a_2=b_2=c_2=d_2=0,$ $d_1=3b_1,$ $b_1=-c_1/4,$  $c_1<0$; 
\item[(b)] $a_1=a_2=b_2=c_1=c_2=d_2=0,$ $d_1=3b_1,$  $b_1<0$;
\item[(c)] $a_2=b_1=b_2=c_2=d_1=d_2=0,$ $a_1^2+c_1<0$;
\item[(d)] $a_1=a_2=b_1=b_2=c_1=c_2=d_1=d_2=0;$
\item[(e)] $a_1=a_2=b_2=c_2=d_2=0,$ $d_1=-b_1-c_1\neq 3b_1,$  $2b_1+c_1\leq 0,$  $b_1>0$;
\item[(f)] $a_2=b_2=c_1=c_2=d_2=0,$ $b_1\neq0,$ $a_1^2+3b_1-d_1<0,$ $a_1^2+4(b_1+d_1)<0;$
\item[(g)] $a_1=a_2=b_1=b_2=c_2=d_2=0,$ $d_1=-c_1>0.$
\end{itemize}
\end{theorem}

The paper is organized as follows. In section \ref{sec:Preliminaries} we present some basic results on the Poincar\'e compactification for polynomial vector fields in the plane $\mathbb{R}^2$, on the characterization of global centers, and on the vertical blow up's. In section \ref{sec:mainresults} we prove our main result (Theorem \ref{main_teo}) about the global centers of system \eqref{main_eq}.

\section{Preliminaries}\label{sec:Preliminaries}

This section is devoted to establishing some basic results that will be used for proving Theorem \ref{main_teo}.

\subsection{Poincaré Compactification for vector fields in the plane}

To study the global dynamics of a planar polynomial differential system $X=(P,Q),$ we must classify the local phase portraits of its finite and infinite equilibrium points in the Poincaré disc.

Let $\sss^2=\{\mathbf{z}\in\mathbb{R}^3:||\mathbf{z}||=1\}$ be the unit sphere of $\mathbb{R}^3.$ We know that the polynomial vector field $X$ induces an analytic vector field in $\sss^2,$ that we denote by $p(X)$, see for instance \cite[Chaper 5]{MR2256001} or \cite{10.2307/2001320}.

The vector field $p(X)$ allows to study the dynamics of the vector field $X$ in the neighbourhood of infinity, that is, in the neighbourhood of the equator $\sss^1=\{\mathbf{z}\in \sss^2:z_3=0\}.$

To obtain the analytical expression for $p(X)$ we shall consider the sphere as a smooth manifold. We choose the six local charts given by $U_i=\{\mathbf{z}\in \sss^2: z_i>0\}$ and $V_i=\{\mathbf{z}\in \sss^2: z_i<0\}$, for $i=1, 2, 3$, with the corresponding coordinate maps $\varphi_i: U_i\rightarrow \mathbb{R}^2$ and $\psi_i: V_i\rightarrow \mathbb{R}^2$, defined by $\varphi_k(\mathbf{z}) = \psi_k(\mathbf{z}) = (z_m/z_k, z_n/z_k)$ for $m < n$ and $m, n\ne k$.
Denote by $(u, v)$ the local coordinates on $U_i$ and $V_i$ for $i=1,2,3.$ From \cite[Chaper 5]{MR2256001} the vector field $p(X)$ in these local charts is
\begin{equation}\label{poincare_comp}
\begin{array}{rl}
(\dot{u},\dot{v})=&\left(v^n\left(-uP\left(\dfrac{1}{v},\dfrac{u}{v}\right)+Q\left(\dfrac{1}{v},\dfrac{u}{v}\right)\right),-v^{n+1}P\left(\dfrac{1}{v},\dfrac{u}{v}\right)\right) \quad\text{in }\quad U_1;\vspace{0.3cm}\\
(\dot{u},\dot{v})=&\left(v^n\left(-uQ\left(\dfrac{u}{v},\dfrac{1}{v}\right)+P\left(\dfrac{u}{v},\dfrac{1}{v}\right)\right),-v^{n+1}Q\left(\dfrac{u}{v},\dfrac{1}{v}\right)\right) \quad\text{in }\quad U_2;\vspace{0.3cm}\\
(\dot{u},\dot{v})=&\left(P(u,v),Q(u,v)\right) \quad\text{in }\quad U_3,
\end{array}
\end{equation}
where $n$ is the degree of the polynomial vector field $X$. Recall that the expressions of the vector field $p(X)$ in the local chart $(V_i,\psi_i)$ is the same that as in the local card $(U_i,\varphi_i)$ multiplied by $(-1)^{n-1}$ for $i=1,2,3.$

The points of the infinity in all the local chart are of the form $(u,0)$. The infinity $\mathbb{S}^1$ is invariant under the flow of $p\left(X\right)$. 		

The equilibrium points of the vector field $X$ are called {\it finite} equilibrium points of $X$ or of $p(X)$, while the equilibrium points of the vector field $p(X)$ in $\mathbb{S}^1$ are called {\it infinite} equilibrium points of $X$ or of $p(X)$.

\subsection{Characterization of the global centers}

The following result gives the necessary and sufficient conditions in order that a planar polynomial differential system has a global center, a proof of this result can be found in \cite{doi:10.1080/14689367.2023.2228737}.

\begin{proposition}\label{prop_main}
A polynomial differential system of degree $n$ in $\mathbb{R}^2$ without a line of equilibrium points at infinity, has a global center if, and only if, it has a unique finite equilibrium point which is a center and all the local phase portraits of the infinite equilibrium points (if they exists) are formed by two hyperbolic sectors, and consequently their two separatrices are on the infinite circle.
\end{proposition}

\subsection{Vertical blow up} 

Consider a real planar polynomial differential system given by
\begin{equation}\label{e1}
\dot{x}=P\left(x, y\right)=P_n\left(x, y\right)+\ldots,\quad
\dot{y}=Q\left(x, y\right)=Q_n\left(x, y\right)+\ldots,
\end{equation}
with $P$ and $Q$ being coprime polynomials, $P_n$ and $Q_n$ being homogeneous polynomials of degree $n\in\mathbb{N}$ and the dots representing higher order terms in $x$ and $y$. Since $n>0$ the origin is an equilibrium point of system \eqref{e1}. Then the \textit{characteristic directions} at the origin are given by the straight lines trought the origin defined by the real linear factors of the homogeneous polynomial 
$$
R_n(x,y)=P_n(x,y)y-Q_n(x,y)x.
$$
It is known that the orbits that start or end at the origin start or end tangent to the straight lines given by the characteristic directions. For more details on the characteristic directions see for instance \cite{MR0350126}.

Suppose that we have an equilibrium point at the origin of coordinates, as in the differential system \eqref{e1} and that this equilibrium is linearly zero. Then for studying its local phase portrait we will do vertical blow up's. 

We define the vertical blow up in the $y$ direction as the change of variables $(u,v)=(x,y/x)$. This change transforms the origin of system \eqref{e1} in the straight line $x=0$, analyzing the dynamics of the differential system in a neighbourhood of this straight line we are analyzing the local phase portrait of the equilibrium point at the origin of system \eqref{e1}. But before doing a vertical blow up in order that we do not lost information we must avoid that the direction $x=0$ be a characteristic direction of the origin of system \eqref{e1}. If $x=0$ is a characteristic direction we do a convenient twist $(x,y)=(u,u+\alpha v)$ with $\alpha\ne 0$.

\section{The classification of global centers}\label{sec:mainresults}

This section is dedicated to prove Theorem \ref{main_teo}. Due to the fact that the proof is very long, we divide it into six parts.

\smallskip

\noindent{\bf Case $(i)$ of Theorem \ref{lemma_main} when $b_1\neq 0$}.
We start by analyzing case $(i)$ of Theorem \ref{lemma_main}. Thus $c_2=0,$ $b_2d_1+d_2b_1=0$ and $d_1=3b_1.$ Since $b_1\neq0$ we get the following system
\begin{equation}\label{main_eq_c1_1_2}
\begin{array}{rl}
\dot{x}=&y+2a_1xy-a_2(x^2-y^2)-\left(b_2-\dfrac{b_2d_1}{b_1}\right)x^3-c_1x^2y+\left(3b_2+\dfrac{b_2d_1}{b_1}\right)xy^2,\\
\dot{y}=&-x+a_1(x^2-y^2)+2a_2xy+(4b_1+c_1)x^3-\left(3b_2+\dfrac{b_2d_1}{b_1}\right)x^2y\\
&+c_1xy^2+\left(b_2-\dfrac{b_2d_1}{b_1}\right)y^3.
\end{array}
\end{equation}
From \eqref{poincare_comp} system \eqref{main_eq_c1_1_2} in the local charts $U_1$ and $U_2$ becomes
\begin{equation}\label{main_eq_c1_1_U1_2}
\begin{array}{rl}
\dot{u}=&4b_1+c_1+2c_1u^2-\dfrac{2b_2d_1u(1+u^2)}{b_1}-2b_2(u+u^3)-v\Big(a_1(-1+3u^2)+v\\
&+u(a_2(-3+u^2)+uv)\Big),\\
\dot{v}=&-\dfrac{v}{b_1} \Big(b_2d_1(1+u^2)+b_1(-c_1u+b_2(-1+3u^2)-a_2v+(2a_1+a_2u+v)uv)\Big),
\end{array}
\end{equation}
and
\begin{equation}\label{main_eq_c1_1_U2_2}
\begin{array}{rl}
\dot{u}=&-2c_1u^2-(4b_1+c_1)u^4+\dfrac{2b_2}{b_1}(b_1+d_1)u(1+u^2)+a_2v-(3a_2u+a_1(-3+u^2))u\\
&v+(1+u^2)v^2,\\
\dot{v}=&v\left(\dfrac{b_2}{b_1}(-b_1+d_1+(3b_1+d_1)u^2)+a_1v\right.-u(c_1+(4b_1+c_1)u^2+(2a_2+a_1u)v-v^2)\Big),
\end{array}
\end{equation}
respectively.

The origin of the chart $U_2$ is an equilibrium point, see system \eqref{main_eq_c1_1_U2_2}. Then the Jacobian matrix of system \eqref{main_eq_c1_1_U2_2} at $(0,0)$ is
$$
\left(\begin{matrix}
\dfrac{2b_2(b_1+d_1)}{b_1} & a_2\\
0 & \dfrac{b_2(d_1-b_1)}{b_1}
\end{matrix}\right).
$$
In order to have a global center of system \eqref{main_eq_c1_1_2} at the origin all the entries of this Jacobian matrix must be zero. This implies that $a_2=b_2=0.$

\smallskip

\begin{figure}[h]
\begin{overpic}[scale=0.45]{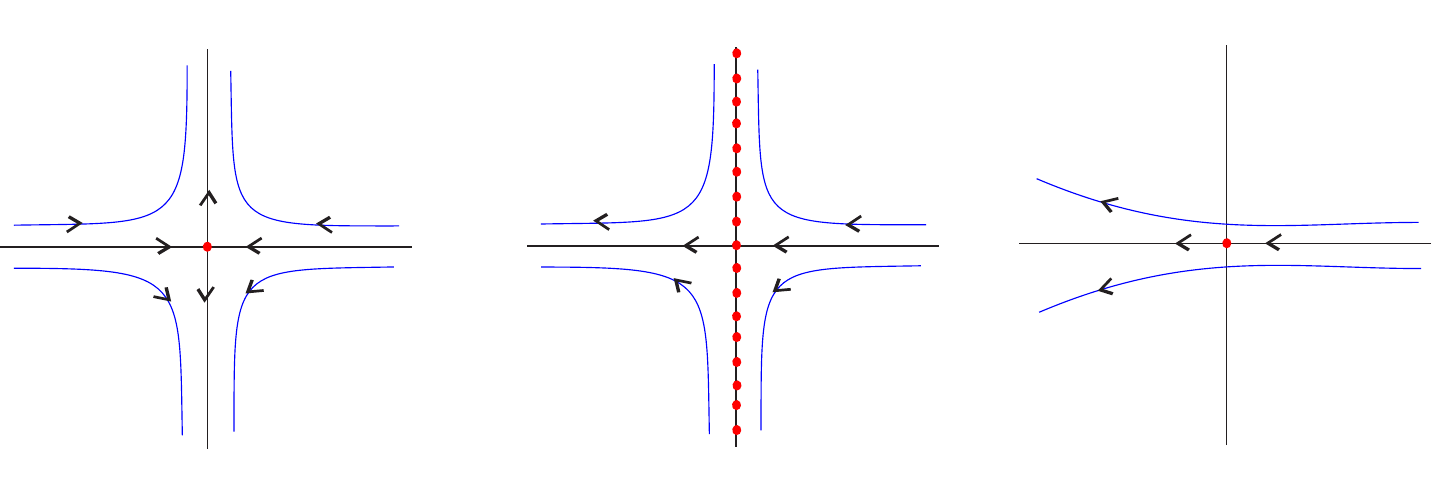}
        \put(101,17.5){$u$}
        \put(66.5,17.5){$u_1$}
        \put(29.5,17.5){$u_1$}
        \put(13.5,32){$v_1$}
        \put(50,32){$v_1$}
        \put(85,32){$v$}
        \put(50,0){$(b)$}
        \put(84,0){$(c)$}
        \put(12.5,0){$(a)$}
\end{overpic}
\caption{\footnotesize{Case (i).1 with $c_1<0$. (a) The local phase portrait at the origin of the differential system \eqref{main_eq_c1_1_U1_blow1_resc_2}. (b) The local phase portrait at the origin of the differential system \eqref{main_eq_c1_1_U1_blow1_2}. (c) The local phase portrait at the origin of the local chart $U_1$ of system \eqref{main_eq_c1_1_2}.}}
\label{fig_4}
\end{figure}
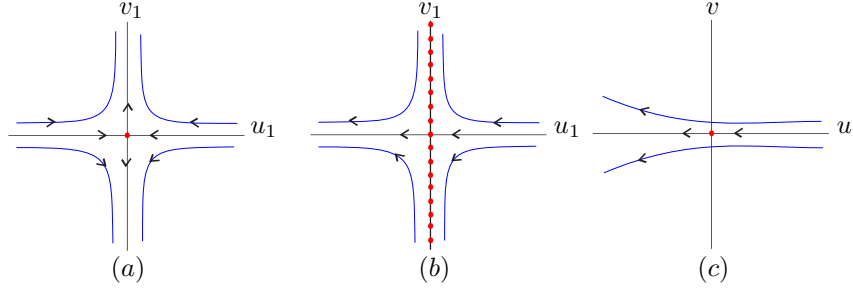

\noindent{\bf Case} (i).1: $c_1\neq 0$. In the chart $U_1$ system \eqref{main_eq_c1_1_U1_2} has the two equilibrium points $p_\pm=\left(\pm\sqrt{(-4b_1-c_1)/(2c_1)},0\right).$ The Jacobian matrix of system \eqref{main_eq_c1_1_U1_2} at $p_\pm$ is
$$
\left(\begin{matrix}
\pm2\sqrt{2}\sqrt{-c_1(4b_1+c_1)}& \dfrac{a_1+3a_1(4b_1+c_1)}{2c_1}\\
0 & \pm\sqrt{\dfrac{-c_1(4b_1+c_1)}{2}}
\end{matrix}\right).
$$
Since the Jacobian matrix at any infinite equilibrium point of system \eqref{main_eq_c1_1_2} must be identically zero if we want a global center, we have that $a_1=0$ and $b_1=-c_1/4.$ This implies that the origin is the unique infinite equilibrium point of the chart $U_1.$ Since $u=0$ is not a characteristic direction at the origin of the chart $U_1$, we do the vertical blow-up $(u,v)=(u_1,u_1v_1).$ Then system \eqref{main_eq_c1_1_U1_2} becomes
\begin{equation}\label{main_eq_c1_1_U1_blow1_2}
\dot{u_1}=u_1^2\Big(2c_1-(1+u_1^2)v_1^2\Big),\quad \dot{v_1}=u_1v_1\Big(v_1^2-c_1\Big).
\end{equation}
Doing a rescaling of the time we eliminate the common factor $u_1$ between $\dot{u_1}$ and $\dot{v_1}$ and we obtain the system
\begin{equation}\label{main_eq_c1_1_U1_blow1_resc_2}
\dot{u_1}=u_1\Big(2c_1-(1+u_1^2)v_1^2\Big),\quad \dot{v_1}=v_1\Big(v_1^2-c_1\Big).
\end{equation}
Then the equilibrium points of system \eqref{main_eq_c1_1_U1_blow1_resc_2} on the straight line $u_1=0$ are
$(0,0)$ and $q_\pm=(0,\pm\sqrt{c_1}).$ The eigenvalues of the Jacobian matrix of system \eqref{main_eq_c1_1_U1_blow1_resc_2} at the origin are $2c_1$ and $-c_1$. If $c_1<0,$ then $(0,0)$ is a hyperbolic saddle and the equilibrium points $q_+$ and $q_-$ do not exist. Therefore going back through the changes of variables we get that the origin of the local chart $U_1$ is formed by two hyperbolic sectors, see Figure \ref{fig_4}.

If $c_1>0$ the eigenvalues of the Jacobian matrix of system \eqref{main_eq_c1_1_U1_blow1_resc_2} at $q_\pm$ are $c_1$ and $2c_1$. Then $q_\pm$ are unstable hyperbolic nodes. Therefore the local phase portraits at the origin of the local chart $U_1$ are not formed by two hyperbolic sectors. Consequently system \eqref{main_eq_c1_1_2} cannot have a global center.

Now we shall study the local phase portrait at the origin of the local chart $U_2$ when $c_1<0$. Since the characteristic directions at the origin of $U_2$ are given by the linear real factors of $v(v^2-c_1u^2)$, we do the vertical blow-up $(u,v)=(u_1,u_1v_1)$ and system \eqref{main_eq_c1_1_U2_2} becomes
\begin{equation*}\label{main_eq_c1_1_U2_blow1_2}
\dot{u_1}=-u_1^2\Big(2c_1-(1+u_1^2)v_1^2\Big),\quad
\dot{v_1}=-u_1v_1\Big(v_1^2-c_1\Big).
\end{equation*}
Doing a rescaling of the time we eliminate the common factor $u_1$ 	 between $\dot{u_1}$ and $\dot{v_1}$ and we obtain the system
\begin{equation}\label{main_eq_c1_1_U2_blow1_resc_2}
\dot{u_1}=-u_1\Big(2c_1-(1+u_1^2)v_1^2\Big),\quad
\dot{v_1}=-v_1\Big(v_1^2-c_1\Big).
\end{equation}
Then the equilibrium points of system \eqref{main_eq_c1_1_U2_blow1_resc_2} on the straight line $u_1=0$ are $(0,0)$ and $q_\pm=(0,\pm\sqrt{c_1})$.
The eigenvalues of the Jacobian matrix of system \eqref{main_eq_c1_1_U2_blow1_resc_2} at the origin are $-2c_1$ and $c_1$. Since $c_1<0,$ then $(0,0)$ is a hyperbolic saddle and the equilibrium points $q_+$ and $q_-$ do not exist. Therefore going back through the changes of variables we get that the local phase portrait at the origin of the local chart $U_2$ is the one of the right picture of Figure \ref{fig_4}, reversing the orientation of the orbits. Hence the origin of the local chart $U_2$ is formed by two hyperbolic sectors. 

In summary, for the values of the parameters $a_1=a_2=b_2=c_2=0,$ $d_1=3b_1,$ $b_1=-c_1/4,$ $b_2d_1+d_2b_1=0$ and $c_1<0$ the infinite equilibrium points satisfy the conditions of Proposition \ref{prop_main} in order to have a global center. Now we must see if the origin of coordinates is the unique finite equilibrium point. 

From the condition $b_2d_1+d_2b_1=0$ and since $b_2=0$ and $b_1\ne 0$ we obtain that $d_2=0$. Therefore the differential system \eqref{main_eq} reduces to
\begin{equation}\label{aa1}
\dot x=y - c_1 x^2 y, \quad \dot y= -x + c_1 x y^2.
\end{equation}
Since the unique real equilibrium of this system is the origin of coordinates by Proposition \ref{prop_main} system \eqref{aa1} has a global center. This completes the proof of statement (a) of Theorem \ref{main_teo}.

\smallskip

\noindent{\bf Case} (i).2: $c_1=0$. Now since in the local chart $U_1$ system \eqref{main_eq_c1_1_U1_2} has no infinite equilibrium points, it is enough to study if the origin of the local chart $U_2$ is an infinite equilibrium point.

Notice that for system \eqref{main_eq_c1_1_U2_2} $R_2(u,v)=v^2(2a_1u+v)$, so $u=0$ is not a characteristic direction. Hence we do the vertical blow-up $(u,v)=(u_1,u_1v_1)$ obtaining the system
\begin{equation}\label{main_eq_c1_1_U2_blow1_resc_trans_blow2_2c1}
\dot{u_1}=-u_1^2\Big(4b_1u_1^2-3a_1v_1+a_1u_1^2v_1-v_1^2-u_1^2v_1^2\Big), \quad \dot{v_1}=-u_1v_1^2\Big(2a_1+v_1\Big).
\end{equation}
Doing a rescaling of the time we eliminate the common factor $u_1$ between $\dot{u_1}$ and $\dot{v_1}$ and we obtain the system
\begin{equation}\label{main_eq_c1_1_U2_blow1_resc_trans_blow2_resc2_2c1}
\dot{u_1}=-u_1\Big(4b_1u_1^2-3a_1v_1+a_1u_1^2v_1-v_1^2-u_1^2v_1^2\Big),\vspace{0.3cm}\\
\dot{v_1}=-v_1^2\Big(2a_1+v_1\Big).
\end{equation}
Then the equilibrium points of system \eqref{main_eq_c1_1_U2_blow1_resc_trans_blow2_resc2_2c1} on the straight line $u_1=0$ are $(0,0)$ and $(0,-2a_1).$ The eigenvalues of the Jacobian matrix of system \eqref{main_eq_c1_1_U2_blow1_resc_trans_blow2_resc2_2c1} at $(0,-2a_1)$ are $-2a_1^2$ and $-4a_1^2$ and the origin is linearly zero.

\noindent\textit{Subcase $1$}: $a_1\neq0.$ Then $(0,-2a_1)$ is a stable hyperbolic node. Therefore the local phase portrait at the origin of the chart $U_2$ is not formed by two hyperbolic sectors, and system \eqref{main_eq_c1_1_2} cannot have a global center.

\noindent\textit{Subcase $2$}: $a_1=0.$ Then the unique equilibrium point of system \eqref{main_eq_c1_1_U2_blow1_resc_trans_blow2_resc2_2c1} on the straight line $u_1=0$ is the $(0,0)$. We analyze its local phase portrait doing blow-up's. For system \eqref{main_eq_c1_1_U2_blow1_resc_trans_blow2_resc2_2c1}
$R_2(u_1,v_1)=5a_1u_1v_1^2$, so $u_1=0$ is a characteristic direction. Consequently before doing a vertical blow-up we translate the direction $u_1=0$ to the direction $u_1=v_1$ doing the change of variables $(u_1,v_1)=(u_2-v_2,v_2).$ In the new variables $(u_2,v_2)$ system \eqref{main_eq_c1_1_U2_blow1_resc_trans_blow2_resc2_2c1} becomes
\begin{equation}\label{main_eq_c1_1_U2_blow1_resc_rot_2c10}
\dot{u_2}=-4b_1(u_2-v_2)^3 +v_2^2(u_2+u_2^3-3u_2^2v_2+3u_2v_2^2-v_2(2+v_2^2)),\quad
\dot{v_2}=-v_2^3.
\end{equation}
Now we do the vertical blow-up $(u_2,v_2)=(u_3,u_3v_3)$ obtaining the system
\begin{equation}\label{main_eq_c1_1_U2_blow1_resc_trans_blow2_2c10}
\begin{array}{l}
\dot{u_3}=-u_3^3\Big(4b_1-12b_1v_3-v_3^2+12b_1v_3^2-u_3^2v_3^2+2v_3^3-4b_1v_3^3+3u_3^2v_3^3-3u_3^2v_3^4+u_3^2v_3^5\Big),\\
\dot{v_3}=u_3^2(v_3-1)v_3\Big(4b_1(1-v_3)^2+(2+u_3^2(-1+v_3)^2)v_3^2\Big).
\end{array}
\end{equation}
Doing a rescaling of the time we eliminate the common factor $u_3^2$ between $\dot{u_3}$ and $\dot{v_3}$ and we obtain the system
\begin{equation}\label{main_eq_c1_1_U2_blow1_resc_trans_blow2_resc2_2c10}
\begin{array}{l}
\dot{u_3}=-u_3\Big(4b_1-12b_1v_3-v_3^2+12b_1v_3^2-u_3^2v_3^2+2v_3^3-4b_1v_3^3+3u_3^2v_3^3-3u_3^2v_3^4+u_3^2v_3^5\Big),\\
\dot{v_3}=(v_3-1)v_3\Big(4b_1(1-v_3)^2+(2+u_3^2(-1+v_3)^2)v_3^2\Big).
\end{array}
\end{equation}
Here we consider two cases.

If $b_1=1/2$, then the equilibrium points of system \eqref{main_eq_c1_1_U2_blow1_resc_trans_blow2_resc2_2c10} on the straight line $u_3=0$ are $(0,0)$, $(0,1)$ and $(0,1/2)$. The eigenvalues of the Jacobian matrix of system \eqref{main_eq_c1_1_U2_blow1_resc_trans_blow2_resc2_2c10} at $(0,1/2)$ are $-1/4$ and $-1,$ then $(0,1/2)$ is a stable hyperbolic node. Therefore the local phase portrait of the origin of the chart $U_2$ is not formed by two hyperbolic sectors. Consequently system \eqref{main_eq_c1_1_2} cannot have a global center.

If $b_1\neq 1/2$, then the equilibrium points of system \eqref{main_eq_c1_1_U2_blow1_resc_trans_blow2_resc2_2c10} on the straight line $u_3=0$ are
$(0,0)$, $(0,1)$ and $q_\pm=(0,(2b_1\pm\sqrt{2b_1})/(2b_1-1))$.

If $b_1>0,$ then the equilibrium points $q_+$ and $q_-$ are  hyperbolic nodes, because
$\lambda_1\lambda_2=(16b_1^2(1+2\sqrt{2b_1}+2b1)^2)/((1-2b_1)^4)>0$,
where $\lambda_1,\lambda_2$ are the eigenvalues of $q_\pm$. Therefore the local phase portrait of the origin of the chart $U_2$ is not formed by two hyperbolic sectors, and system \eqref{main_eq_c1_1_2} cannot have a global center.

\begin{figure}[h]
\begin{overpic}[scale=0.5]{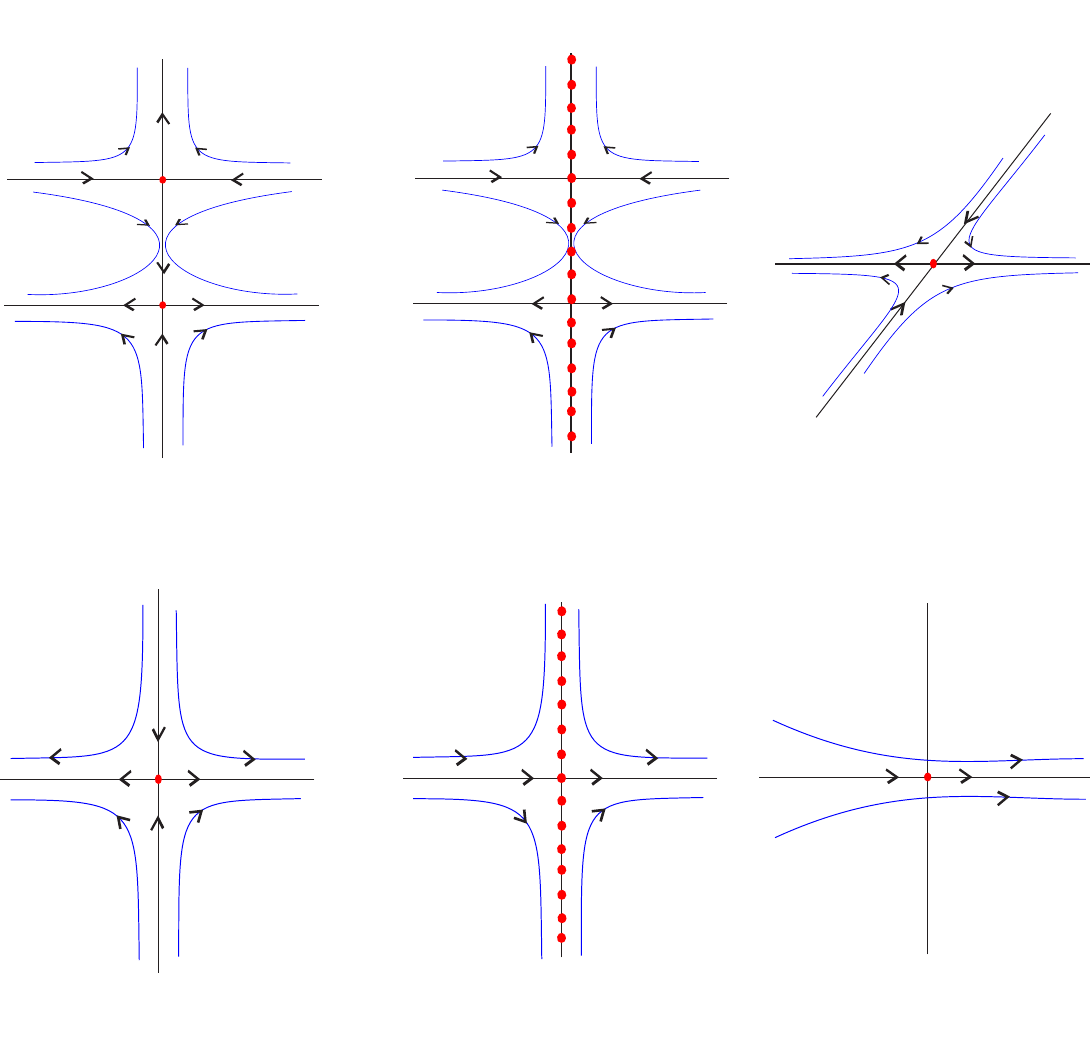}
        \put(101,24){$u$}
        \put(94,88){$v_2=u_2$}
        \put(66,24){$u_1$}
          \put(13.5,43.5){$v_1$}
        \put(67.5,68){$u_3$}
        \put(67.5,79){\tiny{$v_3=1$}}
        \put(50,92.5){$v_3$}
        \put(101,72){$u_2$}
        \put(84.5,43.5){$v$}
        \put(30,24){$u_1$}
        \put(51,43.5){$v_1$}
        \put(30,68){$u_3$}
        \put(30,79){\tiny{$v_3=1$}}
        \put(14,92.5){$v_3$}
          \put(50,50){$(b)$}
        \put(83.5,50){$(c)$}
        \put(12.5,50){$(a)$}
        \put(49.5,2){$(e)$}
        \put(83,2){$(f)$}
        \put(12,2){$(d)$}
		\end{overpic}
\caption{\footnotesize{System \eqref{main_eq_c1_1_2} with $a_1=a_2=b_2=c_1=d_2=0$ and $b_1<0$. (a) The local phase portrait along the straight line $u_3=0$ of the differential system \eqref{main_eq_c1_1_U2_blow1_resc_trans_blow2_resc2_2c10}. (b) The local phase portrait along the straight line $u_3=0$ of the differential system \eqref{main_eq_c1_1_U2_blow1_resc_trans_blow2_2c10} where now the straight line $u_3=0$ is filled with equilibria. (c) The local phase portrait at the origin of the differential system \eqref{main_eq_c1_1_U2_blow1_resc_rot_2c10}. (d) The local phase portrait at the origin of the differential system \eqref{main_eq_c1_1_U2_blow1_resc_trans_blow2_resc2_2c1}. (e) The local phase portrait at the origin of the differential system \eqref{main_eq_c1_1_U2_blow1_resc_trans_blow2_2c1}. (f) The local phase portrait at the origin of the local chart $U_2$ of system \eqref{main_eq_c1_1_U2_2}.}}
\label{fig_6}
\end{figure}
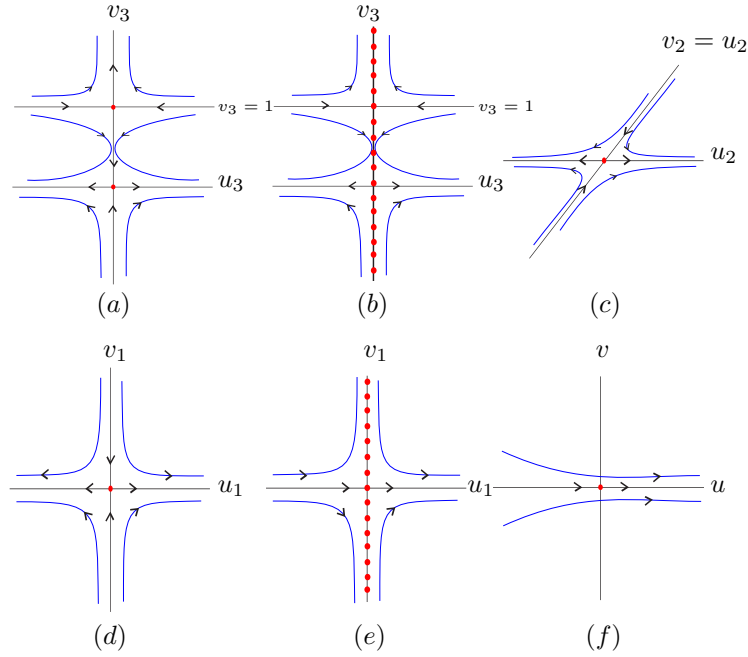

If $b_1<0,$ then the equilibrium points $q_+$ and $q_-$ do not exist, and the equilibrium point $(0,0)$ (resp. $(0,1)$) is a hyperbolic saddle because its associated eigenvalues are $-4b_1$ and $4b_1$ (resp. $-1$ and $2$). Therefore going back through the changes of variables we get that the origin of the local chart $U_2$ is formed by two hyperbolic sectors, see Figure \ref{fig_6}. 

In summary, for the values of the parameters $a_1=a_2=b_2=c_1=c_2=0,$ $d_1=3b_1,$ $b_2d_1+d_2b_1=0$ and $b_1<0$ the infinite equilibrium points satisfy the conditions of Proposition \ref{prop_main} in order to have a global center. Now we must see if the origin of coordinates is the unique finite equilibrium point. 

From the condition $b_2d_1+d_2b_1=0$ and since $b_2=0$ and $b_1\ne 0$ we obtain that $d_2=0$. Therefore the differential system \eqref{main_eq} reduces to
\begin{equation}\label{aa2}
\dot x=y, \quad \dot y= -x + 4 b_1 x^3.
\end{equation}
Since the unique real equilibrium of this system is the origin of coordinates by Proposition \ref{prop_main} system \eqref{aa2} has a global center. This completes the proof of statement (b) of Theorem \ref{main_teo}.

\smallskip

\noindent{\bf Case $(i)$ of Theorem \ref{lemma_main} when $b_1= 0$}.\label{case_1_1}
From Theorem \ref{lemma_main} we know that $c_2=0,$ $b_2d_1+d_2b_1=0,$ $d_1=3b_1.$ Therefore the differential system \eqref{main_eq} reduces to
\begin{equation}\label{main_eq_c1_1}
\begin{array}{l}
\dot{x}=y+2a_1xy-a_2(x^2-y^2)-(b_2+d_2)x^3-c_1x^2y+(3b_2-d_2)xy^2,\\
\dot{y}=-x+a_1(x^2-y^2)+2a_2xy+c_1x^3-(3b_2-d_2)x^2y+c_1xy^2+(b_2+d_2)y^3.
\end{array}
\end{equation}
From \eqref{poincare_comp} system \eqref{main_eq_c1_1} in the local charts $U_1$ and $U_2$ becomes
\begin{equation}\label{main_eq_c1_1_U1}
\begin{array}{l}
\dot{u}=c_1+2c_1u^2-2(b_2-d_2)(u+u^3)+a_1v-(3a_1u+a_2(-3+u^2))uv-(1+u^2)v^2,\\
\dot{v}=v \Big(b_2+d_2+c_1u-3b_2u^2+a_2v+u(d_2u-v(2a_1+a_2u+v))\Big),
\end{array}
\end{equation}
and
\begin{equation}\label{main_eq_c1_1_U2}
\begin{array}{rl}
\dot{u}=&-c_1u^2(2+u^2)+2b_2(u+u^3)-2d_2(u+u^3)+a_2v\\
&-(3a_2u+a_1(-3+u^2))uv+(1+u^2)v^2,\\
\dot{v}=&-v\Big(b_2-3b_2u^2+(d_2+c_1u)(1+u^2)+2a_2uv+a_1(-1+u^2)v-uv^2\Big),
\end{array}
\end{equation}
respectively.

In the chart $U_2$ system \eqref{main_eq_c1_1_U2} has the origin as an infinite equilibrium point. Then the Jacobian matrix of system \eqref{main_eq_c1_1_U2} at $(0,0)$ is
$$
\left(\begin{matrix}
2b_2-2d_2 & a_2\\
0 & -b_2-d_2
\end{matrix}\right).
$$
In order to have a global center of system \eqref{main_eq_c1_1} at the origin of coordinates, all the entries of this Jacobian matrix must be zero. This implies that $a_2=b_2=d_2=0.$ With these new conditions we have that the origin of the chart $U_2$ is the unique infinite equilibrium point of system \eqref{main_eq_c1_1}, because in the chart $U_1$ system \eqref{main_eq_c1_1_U1} has no equilibrium points. Hence it is enough to study the local phase portrait at the origin of the local chart $U_2.$

We need to do blow up's for studying the local phase portrait at the origin of the chart $U_2$. Since $R_2(u,v)=v(-c_1u^2+2a_1uv+v^2)$ the straight line $u=0$ is not a characteristic direction and we do the vertical blow-up $(u,v)=(u_1,u_1v_1).$ Then system \eqref{main_eq_c1_1_U2} becomes
\begin{equation}\label{main_eq_c1_1_U2_blow1}
\begin{array}{l}
\dot{u_1}=u_1^2\Big(-c_1(2+u_1^2)+v_1(3a_1-a_1u_1^2+v_1+u_1^2v_1)\Big),\\
\dot{v_1}=u_1v_1\Big(c_1-v_1(2a_1+v_1)\Big).
\end{array}
\end{equation}
Doing a rescaling of the time we eliminate the common factor $u_1$ between $\dot{u_1}$ and $\dot{v_1}$ and we obtain the system
\begin{equation}\label{main_eq_c1_1_U2_blow1_resc}
\begin{array}{l}
\dot{u_1}=u_1\Big(-c_1(2+u_1^2)+v_1(3a_1-a_1u_1^2+v_1+u_1^2v_1)\Big),\\
\dot{v_1}=v_1\Big(c_1-v_1(2a_1+v_1)\Big).
\end{array}
\end{equation}
Then the equilibrium points of system \eqref{main_eq_c1_1_U2_blow1_resc} on the straight line $u_1=0$ are $(0,0)$ and $q_\pm=\left(0,-a_1\pm\sqrt{a_1^2+c_1}\right)$. The eigenvalues of the Jacobian matrix of system \eqref{main_eq_c1_1_U2_blow1_resc} at the origin are $-2c_1$ and $c_1$. We consider three subcases.

\noindent\textit{Subcase $1$}: $a_1^2+c_1<0.$ Then $c_1<0$ and $(0,0)$ is a hyperbolic saddle and the equilibrium points $q_+$ and $q_-$ do not exist. Therefore going back through the changes of variables we get that the origin of the local chart $U_2$ is formed by two hyperbolic sectors, see Figure \ref{fig_4} reversing the sense of all the orbits. 

In summary, for the values of the parameters $a_2=b_1=b_2=c_2=d_1=d_2=0$ and $a_1^2+c_1<0$ the infinite equilibrium points satisfy the conditions of Proposition \ref{prop_main} in order to have a global center. Now we must see if the origin of coordinates is the unique finite equilibrium point. Under the previous conditions the differential system \eqref{main_eq} reduces to
\begin{equation}\label{aa3}
\dot x=y + 2 a_1 x y - c_1 x^2 y,\quad \dot y= -x+ a_1 (x^2 - y^2) + c_1 x^3 + c_1 x y^2. 
\end{equation}
Since the unique real equilibrium of this system is the origin of coordinates by Proposition \ref{prop_main} system \eqref{aa3} has a global center. This completes the proof of statement (c) of Theorem \ref{main_teo}.

\smallskip

\noindent\textit{Subcase $2$}: $a_1^2+c_1>0.$ Then the eigenvalues of the Jacobian matrix of system \eqref{main_eq_c1_1_U2_blow1_resc} at $q_\pm$ are $\lambda_1^\pm=-(a_1^2+c_1\mp a_1\sqrt{a_1^2+c_1})$ and $\lambda_2^\pm=-2(a_1^2+c_1\mp a_1\sqrt{a_1^2+c_1})$.

If $c_1\neq 0$, then $\lambda_1^\pm\lambda_2^\pm=2\left(a_1^2+c_1\mp a_1\sqrt{a_1^2+c_1}\right)^2>0,$ so $q_+$ and $q_-$ are hyperbolic nodes. Therefore the origin of the local chart $U_2$ is not formed by two hyperbolic sectors. Consequently system \eqref{main_eq_c1_1} cannot have a global center.

If $c_1=0$, then the equilibrium points are $(0,0)$
and $(0,-2a_1)$, and the eigenvalues of the Jacobian matrix of system \eqref{main_eq_c1_1_U2_blow1_resc} at $(0,-2a_1)$ are $-2a_1^2$ and $-4a_1^2,$ which implies that $(0,-2a_1)$ is a stable hyperbolic node. Therefore system \eqref{main_eq_c1_1} cannot have a global center.

\noindent\textit{Subcase $3$}: $a_1^2+c_1=0.$ Then the equilibrium points $q_+=q_-=(0,-a_1).$

\begin{figure}[h]
	\begin{overpic}[scale=0.45]{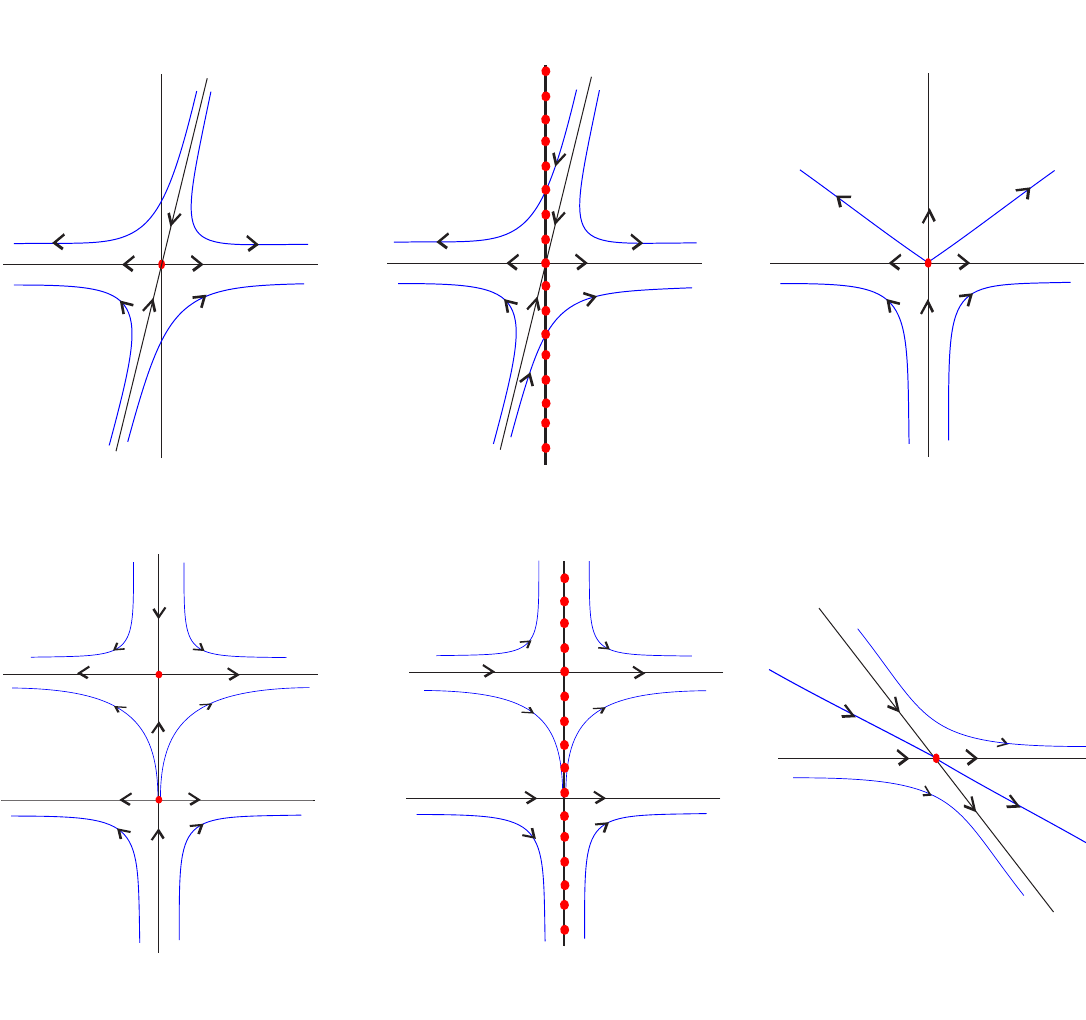}
		\put(101,23.5){$u$}
		\put(94,4){$v=-a_1u$}
		\put(66.5,31.5){$u_1$}
		\put(66,20){$-a_1$}
		\put(13.5,43.5){$v_1$}
		\put(65.5,70){$u_3$}
		\put(49,89){$v_3$}
		\put(101,70){$u_2$}
		\put(84.5,89){$v_2$}
		\put(29,20){$-a_1$}
		\put(30,31){$u_1$}
		\put(51,43.5){$v_1$}
		\put(30,70){$u_3$}
		\put(14,89){$v_3$}
		\put(50,48){$(b)$}
		\put(83.5,48){$(c)$}
		\put(12.5,48){$(a)$}
		\put(50,0){$(e)$}
		\put(85,0){$(f)$}
		\put(12,0){$(d)$}
	\end{overpic}
	\caption{\footnotesize{System \eqref{main_eq_c1_1} with $a_2=b_2=d_2=0$ and $a_1>0$. (a) The local phase portrait at the origin of system \eqref{main_eq_c1_1_U2_blow1_resc_trans_blow2_resc2}. (b) The local phase portrait at the origin of system \eqref{main_eq_c1_1_U2_blow1_resc_trans_blow2}. (c) The local phase portrait at the origin of system \eqref{main_eq_c1_1_U2_blow1_resc_trans}. (d) The local phase portrait in a neighbourhood of the straight line $u_1=0$ of system \eqref{main_eq_c1_1_U2_blow1_resc}. (e) The local phase portrait in a neighbourhood of the straight line $u_1=0$ of system \eqref{main_eq_c1_1_U2_blow1}. (f) The local phase portrait at the origin of the local chart $U_2$ of system \eqref{main_eq_c1_1_U2}.}}
	\label{fig_2}
\end{figure}
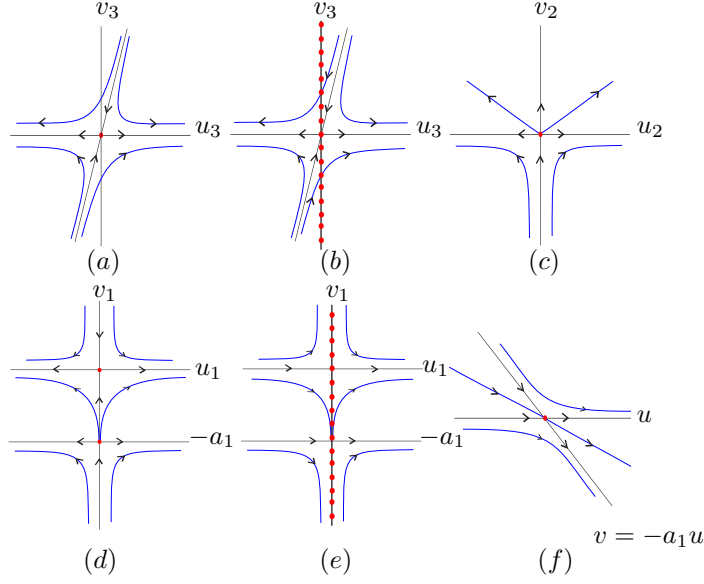

If $a_1\neq 0,$ then the equilibrium point $(0,0)$ is a hyperbolic saddle and $q_+$ is linearly zero. Doing the change of variables $(u_1,v_1)=(u_2,v_2-a_1)$, we translate the  equilibrium point $(0,-a_1)$ to the origin of coordinates. Then system \eqref{main_eq_c1_1_U2_blow1_resc} in the variables $(u_2,v_2)$ becomes
\begin{equation}\label{main_eq_c1_1_U2_blow1_resc_trans}
\begin{array}{l}
\dot{u_2}=u_2\Big(3a_1^2u_2^2+(1+u_2^2)v_2^2+a_1(v_2-3u_2^2v_2)\Big),\\
\dot{v_2}=v_2^2\Big(a_1-v_2\Big).
\end{array}
\end{equation}
Since $R_2(u,v)\equiv0$ all the directions are characteristic, in this case we do the blow-up $(u_2,v_2)=(u_3,u_3v_3)$ and system \eqref{main_eq_c1_1_U2_blow1_resc} writes
\begin{equation}\label{main_eq_c1_1_U2_blow1_resc_trans_blow2}
\begin{array}{l}
\dot{u_3}=u_3^2\Big(3a_1^2u_3+(1+u_3^2)u_3v_3^2+a_1(v_3-3u_3^2v_3)\Big),\\
\dot{v_3}=-u_3^2v_3\Big(3a_1^2-3a_1u_3v_3+(2+u_3^2)v_3^2\Big).
\end{array}
\end{equation}
Doing a rescaling of the time we eliminate the common factor $u_3^2$ between $\dot{u_3}$ and $\dot{v_3}$ and we obtain the system
\begin{equation}\label{main_eq_c1_1_U2_blow1_resc_trans_blow2_resc2}
\begin{array}{l}
\dot{u_3}=3a_1^2u_3+(1+u_3^2)u_3v_3^2+a_1(v_3-3u_3^2v_3),\\
\dot{v_3}=-v_3\Big(3a_1^2-3a_1u_3v_3+(2+u_3^2)v_3^2\Big).
\end{array}
\end{equation}
Then the unique equilibrium point of system \eqref{main_eq_c1_1_U2_blow1_resc_trans_blow2_resc2} on the straight line $u_3=0$ is the $(0,0).$ The eigenvalues of the Jacobian matrix of system \eqref{main_eq_c1_1_U2_blow1_resc_trans_blow2_resc2} at the origin are $3a_1^2$ and $-3a_1^2$. Hence the equilibrium point $(0,0)$ is a hyperbolic saddle. Therefore going back through the changes of variables when $a_1>0$ (the case $a_1<0$ is similar) the local phase portrait at the origin of the local chart $U_2$ is shown in Figure \ref{fig_2}. Consequently there are orbits which go or come from infinity in system \eqref{main_eq_c1_1}, and thus the center of this system cannot global.

If $a_1=0,$ then $c_1=0$ and all the parameters are 0 and we are in the trivial case in which system \eqref{main_eq_c1_1} is $\dot{x}=y$, $\dot{y}=-x.$ Thus we have that the origin is a global center. This completes the proof of statement (d) of Theorem \ref{main_teo}.

\smallskip

\noindent{\bf Case $(ii)$ of Theorem \ref{lemma_main} when $b_1\neq 0$}.
We have that $c_2=0,$ $b_2d_1+d_2b_1=0,$ $d_1\neq 3b_1$ and $F=a_2^2 d_2^3 - 3 a_2^2 d_2d_1^2 + 6 a_2a_1d_2^2 d_1 - 2 a_2a_1d_1^3 - a_1^2 d_2^3 + 3 a_1^2 d_2d_1^2=0.$ Since $b_1\neq0$ system \eqref{main_eq} reduces to system
\begin{equation}\label{main_eq_c22_1}
\begin{array}{rl}
\dot{x}=&y+2a_1xy-\left(b_2-\dfrac{b_2d_1}{b_1}\right)x^3-(3b_1+c_1-d_1)x^2y+\left(3b_2+\dfrac{b_2d_1}{b_1}\right)xy^2\\
&-a_2(x^2-y^2),\\
\dot{y}=&-x+(b_1+c_1+d_1)x^3+2a_2xy-\left(3b_2+\dfrac{b_2d_1}{b_1}\right)x^2y+(-3b_1+c_1+d_1)xy^2\\
&+\left(b_2-\dfrac{b_2d_1}{b_1}\right)y^3+a_1(x^2-y^2).\\
\end{array} 
\end{equation}
From \eqref{poincare_comp}, system \eqref{main_eq_c22_1} in the local charts $U_1$ and $U_2$ becomes
 \begin{equation}\label{main_eq_c22_1_U1}
\begin{array}{rl}
\dot{u}=&b_1+c_1+d_1+2c_1u^2-\dfrac{2b_2d_1u(1+u^2)}{b_1}-2b_2(u+u^3)+a_1v\\
&-v\Big(a_2u(-3+u^2)+v+u^2(3a_1+v)\Big),\\
\dot{v}=&v\Big(b_2+3b_1u-3b_2u^2-\dfrac{b_2d_1(1+u^2)}{b_1}+a_2v
\\
&-u(-c_1+d_1+v(2a_1+a_2u+v))\Big),\\
\end{array}
\end{equation}
and 
 \begin{equation}\label{main_eq_c22_1_U2}
\begin{array}{rl}
\dot{u}=&-2c_1u^2-(b_1+c_1+d_1)u^4+\dfrac{2b_2}{b_1}(b_1+d_1)u(1+u^2) +a_2v\\
&-u(3a_2u+a_1(-3+u^2))v+(1+u^2)v^2,\\
\dot{v}=&v \Big(\dfrac{b_2}{b_1}(-b_1+d_1+(3b_1+d_1)u^2)+a_1v-u(c_1+d_1+(c_1+d_1)u^2\\
&+b_1(-3+u^2)+2a_2v+a_1uv-v^2)\Big),\\
\end{array}
\end{equation}
respectively.

In the chart $U_2$ system \eqref{main_eq_c22_1_U2} has the origin as an equilibrium point. The Jacobian matrix of system \eqref{main_eq_c22_1_U2} at $(0,0)$ is 
$$\left(\begin{matrix}
\dfrac{2b_2(d_1+b_1)}{b_1} & a_2\\
0 & \dfrac{b_2(d_1-b_1)}{b_1}
\end{matrix}\right).$$
In order to have a global center of system \eqref{main_eq_c22_1} at the origin, all entries of this Jacobian matrix must be 0. This implies that $a_2=b_2=0.$ Note that these new conditions implies that $F=0.$

\smallskip

\noindent{\bf Case} (ii).1:\label{sec:case22}
$c_1\neq 0$. Then in chart $U_1$ system \eqref{main_eq_c22_1_U1} has two equilibrium points: $p_\pm=\left(\pm\sqrt{(-b_1-c_1-d_1)/(2c_1)},0\right).$ The Jacobian matrix of system \eqref{main_eq_c22_1_U1} at $p_\pm$ is
$$
\left(\begin{matrix}
\pm2\sqrt{2}\sqrt{-c_1(b_1+c_1+d_1)}& \dfrac{a_1(3b_1+5c_1+3d_1)}{2c_1}\\
0 & \pm(3b_1+c_1-d_1)\sqrt{\dfrac{-(b_1+c_1+d_1)}{2c_1}}
\end{matrix}\right).
$$
Since the Jacobian matrix at any infinite equilibrium point of system \eqref{main_eq_c22_1} must be identically zero, we obtain that $a_1=0$ and $d_1=-b_1-c_1.$ This implies that the origin is the unique equilibrium point in the chart $U_1.$ To determine the local phase portrait at the origin we must do blow up's. The straight line $u=0$ is not a characteristic direction at the origin of $U_1$ because $R_2(u,v)=-v(v^2+4b_1u^2).$ Hence we can do the vertical blow up $(u,v)=(u_1,u_1v_1)$, and system \eqref{main_eq_c22_1_U1} becomes 
\begin{equation}\label{main_eq_c222_1_U1_blow1_3}
\dot{u_1}=u_1^2\Big(2c_1-(1+u_1^2)v_1^2\Big),\quad \dot{v_1}=u_1\Big(4b_1v_1+v_1^3\Big).
\end{equation}
Doing a rescaling of the time we eliminate the common factor $u_1$ between $\dot{u_1}$ and $\dot{v_1}$ and we obtain the system
\begin{equation}\label{main_eq_c222_1_U1_blow1_resc_3}
\dot{u_1}=u_1\Big(2c_1-(1+u_1^2)v_1^2\Big),\quad \dot{v_1}=4b_1v_1+v_1^3.
\end{equation}
The equilibrium points of this system on the straight line $u_1=0$ are $(0,0)$ and $q_\pm=(0,\pm2\sqrt{-b_1}).$
The eigenvalues of the Jacobian matrix of system \eqref{main_eq_c222_1_U1_blow1_resc_3} at the origin are $2c_1$ and $4b_1$. 

Assume $b_1>0$. If $c_1>0$ the origin of system \eqref{main_eq_c222_1_U1_blow1_resc_3} is a node, and consequently going back with the blow down the origin of $U_1$ cannot have two hyperbolic sectors, and we cannot have a global center. Suppose that  $c_1<0$. Then the origin is the unique equilibrium point, which is a hyperbolic saddle. Therefore, going back through the changes of variables we get that the origin of the local chart $U_1$ is formed by two hyperbolic sectors. The pictures of the blow down
in this case coincide with the ones of Figure \ref{fig_4}.

Assume now that $b_1<0$. If $c_1<0$ the origin is a node, and as in the case $b_1>0$ we cannot have a global center. So consider $c_1>0$. Then the eigenvalues of the Jacobian matrix of system \eqref{main_eq_c222_1_U1_blow1_resc_3} at $q_\pm$ are $4b_1+2c_1$ and $-8b_1.$

\begin{figure}[h]
\begin{overpic}[scale=0.55]{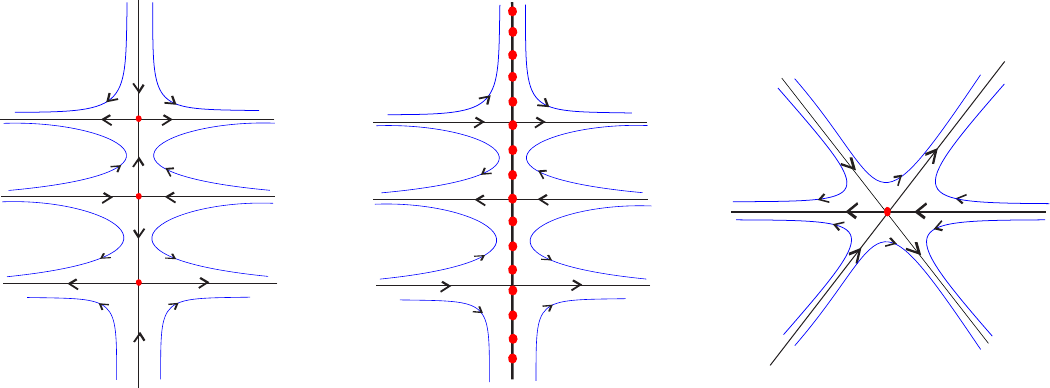}
        \put(101,16){$u$}
        \put(62,18){$u_1$}
        \put(27,18){$u_1$}
        \put(12.5,37.5){$v_1$}
        \put(48,37.5){$v_1$}
        \put(92,32){\tiny{$v=2\sqrt{-b_1}u$}}
        \put(69,32){\tiny{$v=-2\sqrt{-b_1}u$}}
        \put(47.5,-2){$(b)$}
        \put(82.5,-2){$(c)$}
        \put(11.5,-2){$(a)$}
		\end{overpic}
\caption{\footnotesize{Case (ii).1 with $c_1\leq-2b_1.$ (a) The local phase portrait at the origin of the differential system \eqref{main_eq_c222_1_U1_blow1_3}. (b) The local phase portrait at the origin of the differential system \eqref{main_eq_c222_1_U1_blow1_resc_3}. (c) The local phase portrait at the origin of the local chart $U_2$ of system \eqref{main_eq_c22_1_U2}.}}
\label{fig_16}
\end{figure}
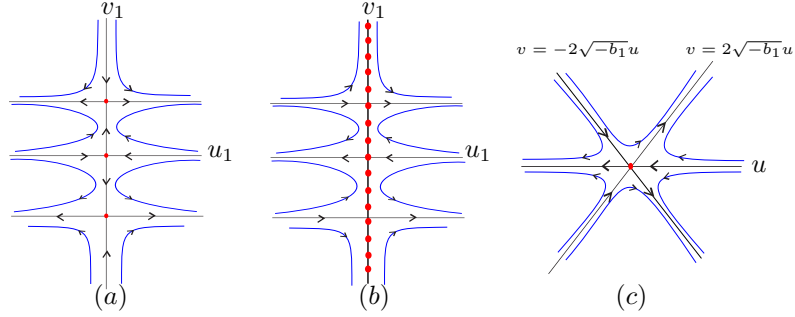

If $c_1=-2b_1$ (resp. $c_1<-2b_1$), then from \cite[Thorem 2.19]{MR2256001} (resp. \cite[Thorem 2.15]{MR2256001}), we can conclude that $q_\pm$ are semi-hyperbolic saddles (resp. hyperbolic saddles). Going back through the changes of variables the local phase portrait at the origin of the chart $U_1$ is shown in Figure \ref{fig_16}. Consequently there are orbits of system \eqref{main_eq_c22_1} which go or come from the infinity, thus the center of this system is not global.
    
If $c_1>-2b_1$ then $q_+$ and $q_-$ are unstable hyperbolic nodes. Consequently, going back through the changes of variables we obtain that there are orbits of system \eqref{main_eq_c22_1} which go or come from the infinity, thus the center of this system cannot global.

In summary, the only possibility of having a global center in this case is when $b_1>0$ and $c_1<0.$ Under this conditions we must study the local phase portrait at the origin of the chart $U_2$ that it is a linearly zero equilibrium point. Since the characteristic directions at the origin of $U_2$ are given by $R_2(u,v)=v(v^2-2(c_1+2b_1)u^2)$ the straight line $u=0$ is not a characteristic direction. Hence we do the vertical blow-up $(u,v)=(u_1,u_1v_1)$, and system \eqref{main_eq_c22_1_U2} becomes 
\begin{equation*}\label{main_eq_c222_1_U2_blow1_31}
\dot{u_1}=u_1^2\Big(-2c_1+(1+u_1^2)v_1^2\Big),\quad \dot{v_1}=u_1v_1\Big(4b_1+2c_1-v_1^2\Big).
\end{equation*}
Doing a rescaling of the time we eliminate the common factor $u_1$ between $\dot{u_1}$ and $\dot{v_1}$ and we obtain the system
\begin{equation}\label{main_eq_c222_1_U2_blow1_resc_31}
\dot{u_1}=u_1\Big(-2c_1+(1+u_1^2)v_1^2\Big),\quad \dot{v_1}=v_1\Big(4b_1+2c_1-v_1^2\Big).
\end{equation}
The equilibrium points of this system on the straight line $u_1=0$ are $(0,0)$ and $q_\pm=\left(0, \pm\sqrt{4b_1+2c_1}\right).$ The eigenvalues of the Jacobian matrix of system \eqref{main_eq_c222_1_U2_blow1_resc_31} at the origin are $-2c_1$ and $4b_1+2c_1$. 

If $c_1=-2b_1$ (resp. $c_1<-2b_1$), then from \cite[Thorem 2.19]{MR2256001} (resp. \cite[Thorem 2.15]{MR2256001}), we can conclude that the origin is a semi-hyperbolic saddle (resp. hyperbolic saddle), its local phase portrait is the one of the left picture of Figure \ref{fig_4}  but with the orbits travelled in converse sense. Going back through the changes of variables the local phase portrait at the origin of system \eqref{main_eq_c222_1_U2_blow1_resc_31} is the one of the right picture of Figure \ref{fig_4} but with the orbit run in reverse sense. Consequently the origin of the local chart $U_2$ is formed by two hyperbolic sectors. 

In summary, for the values of the parameters $a_1=a_2=b_2=c_2=0,$ $d_1=-b_1-c_1\neq 3b_1,$  $2b_1+c_1\leq 0,$ $b_2d_1+d_2b_1=0$ and $b_1>0$ the infinite equilibrium points satisfy the conditions of Proposition \ref{prop_main} in order to have a global center. Now we must see if the origin of coordinates is the unique finite equilibrium point. 

From the condition $b_2d_1+d_2b_1=0$ and since $b_2=0$ and $b_1> 0$ we obtain that $d_2=0$. Therefore the differential system \eqref{main_eq} reduces to
\begin{equation}\label{aa4}
\dot x=y -(4 b_1 + 2 c_1) x^2 y,\quad \dot y= -x - 4 b_1 x y^2.
\end{equation}
Since the unique real equilibrium of this system is the origin of coordinates by Proposition \ref{prop_main} system \eqref{aa4} has a global center. This completes the proof of statement (e) of Theorem \ref{main_teo}.

If $c_1>-2b_1$ then the origin is an unstable hyperbolic node. Consequently going back through the changes of variables we obtain that there are orbits of system \eqref{main_eq_c22_1} which go or come from the infinity, thus the center of this system cannot global.

\smallskip

\noindent{\bf Case} (ii).2:\label{case2.221} 
$c_1=0$ and $b_1+d_1\neq 0$. Then in the chart $U_1$ system \eqref{main_eq_c22_1_U1} has no infinite equilibrium points. So we study the equilibrium point at the origin of the chart $U_2.$ 

Since for system \eqref{main_eq_c22_1_U2} it follows that $R_2(u,v)=v^2(2a_1u+v+(-3b_1+d_1)u^2)$ the straight line $u=0$ is not a characteristic direction. Doing the vertical blow-up $(u,v)=(u_1,u_1v_1)$ system \eqref{main_eq_c22_1_U2} becomes
\begin{equation*}\label{main_eq_c22_1_U2_blow1_resc_trans_blow2_2c1}
\begin{array}{rl}
\dot{u_1}=&u_1^2\Big(-(b_1+d_1)u_1^2-a_1(-3+u_1^2)v_1+(1+u_1^2)v_1^2\Big),\\
\dot{v_1}=&-u_1v_1\Big(-3b_1+d_1+v_1(2a_1+v_1)\Big).
\end{array} 
\end{equation*}
Doing a rescaling of the time we eliminate the common factor $u_1$ between $\dot{u_1}$ and $\dot{v_1}$ and we obtain the system
\begin{equation}\label{main_eq_c22_1_U2_blow1_resc_trans_blow2_resc2_2c1a}
\begin{aligned}
\begin{array}{rcl}
\dot{u_1}&=&u_1\Big(-(b_1+d_1)u_1^2-a_1(-3+u_1^2)v_1+(1+u_1^2)v_1^2\Big),\vspace{0.3cm}\\
\dot{v_1}&=&-v_1\Big(-3b_1+d_1+v_1(2a_1+v_1)\Big).
\end{array} 
\end{aligned}
\end{equation}
The equilibrium points of system \eqref{main_eq_c22_1_U2_blow1_resc_trans_blow2_resc2_2c1a} on the straight line $u_1=0$ are $(0,0)$ and $q_\pm=\Big(0,-a_1\pm\sqrt{a_1^2+3b_1-d_1}\Big).$ The eigenvalues of the Jacobian matrix of system \eqref{main_eq_c22_1_U2_blow1_resc_trans_blow2_resc2_2c1a} at the origin are $0$ and $3b_1-d_1$. 

\noindent\textit{Subcase $1$}: $a_1^2+3b_1-d_1<0.$ Then $3b_1-d_1<0$ and the equilibrium points $q_+$ and $q_-$ do not exist. If $b_1+d_1<0$ (resp. $b_1+d_1>0$), then \cite[Thorem 2.19]{MR2256001} implies that the origin is a semi-hyperbolic saddle (resp. a semi-hyperbolic node as we know this case cannot produce a global center). Going back through the changes of variables we get that the origin of the local chart $U_2$ is formed by two hyperbolic sectors when $b_1+d_1<0$. The pictures of the blow down in this case coincide with the ones of Figure \ref{fig_4} reversing the sense of the orbits. 

In summary, for the values of the parameters $a_2=b_2=c_1=c_2=0,$ $b_1\neq0,$ $b_1+d_1\neq 0,$ $b_2d_1+d_2b_1=0$ and $a_1^2+3b_1-d_1<0,$ the infinite equilibrium points satisfy the conditions of Proposition \ref{prop_main} in order to have a global center. Now we must see if the origin of coordinates is the unique finite equilibrium point. 

From the condition $b_2d_1+d_2b_1=0$ and since $b_2=0$ and $b_1\ne 0$ we obtain that $d_2=0$. Therefore the differential system \eqref{main_eq} reduces to
\begin{equation}\label{bb5}
\dot x=y+2a_1 xy-(3b_1-d_1)x^2y,\quad \dot y=-x+a_1(x^2-y^2)+(b_1 + d_1)x^3+(d_1-3 b_1)xy^2.
\end{equation}
The unique equilibrium of this system is the $(0,0)$ if $a_1^2+4(b_1+d_1)<0$. Therefore by Proposition \ref{prop_main} system \eqref{bb5} has a global center. This completes the proof of statement (f) of Theorem \ref{main_teo}.

\noindent\textit{Subcase $2$}: $a_1^2+3b_1-d_1=0.$ Then the equilibrium points $q_+=q_-=(0,-a_1).$ The eigenvalues of the Jacobian matrix of system \eqref{main_eq_c22_1_U2_blow1_resc_trans_blow2_resc2_2c1a} at $q_+$ are $-a_1^2$ and 0. Recall that $a_1\neq0$ because $d_1\neq 3b_1.$ Thus $3b_1-d_1<0$ and by \cite[Thorem 2.19]{MR2256001}, we can conclude that $q_+$ is a semi-hyperbolic saddle-node. Therefore there are orbits of system \eqref{main_eq_c22_1} which go or come from the infinity. Consequently system \eqref{main_eq_c22_1} cannot have a global center.

\noindent\textit{Subcase $3$}: $a_1^2+3b_1-d_1>0.$ Then system \eqref{main_eq_c22_1_U2_blow1_resc_trans_blow2_resc2_2c1a} has three equilibrium points on $u_1=0$ and doing an analysis similar to case (ii).1 when $c_1\leq-2b_1$ there are orbits of system \eqref{main_eq_c22_1} which go or come from the infinity, thus the center of this system cannot global, see Figure \ref{fig_16}.

\noindent{\bf Case} (ii).3: $c_1=0$ and $b_1+d_1=0.$ From systems  \eqref{main_eq_c22_1_U1} and \eqref{main_eq_c22_1_U2} the circle of the infinity is filled with equilibria. Doing a rescaling of the time we eliminate the common factor $v$ between $\dot{u}$ and $\dot{v}$ in the chart $U_1$ and we obtain the system
\begin{equation}\label{main_eq_c282_1_U1_resc}
\dot{u}=a_1(1-3u^2)-v-u^2v,\quad
\dot{v}=4b_1u-uv(2a_1+v).
\end{equation}
 Since $b_1\neq0,$ then $\dot{v}\Big|_{v=0}=4b_1u\neq 0.$ This means that depending on the sign of $b_1$, there are orbits which go or come from infinity in system \eqref{main_eq_c282_1_U1_resc}, and the center cannot global. 

\smallskip

\noindent{\bf Case $(2)$ of Theorem \ref{lemma_main} when $b_1=0$}.\label{case_2_1}
From Theorem \ref{lemma_main} we know that $c_2=0,$ $b_2d_1+d_2b_1=0,$ $d_1\neq3b_1$ and $F=a_2^2 d_2^3 - 3 a_2^2 d_2d_1^2 + 6 a_2a_1d_2^2 d_1 - 2 a_2a_1d_1^3 - a_1^2 d_2^3 + 3 a_1^2 d_2d_1^2=0.$ Using that $c_2=b_1=0$ system \eqref{main_eq} becomes
\begin{equation}\label{main_eq_c2_1}
\begin{array}{l}
\dot{x}=-(b_2+d_2)x^3+y+2a_1xy-(c_1-d_1)x^2y+(3b_2-d_2)xy^2-a_2(x^2-y^2),\\
\dot{y}=-x+(c_1+d_1)x^3+2a_2xy-(3b_2-d_2)x^2y+(c_1+d_1)xy^2+(b_2+d_2)y^3+a_1(x^2-y^2).
\end{array} 
\end{equation}
From \eqref{poincare_comp} system \eqref{main_eq_c2_1} in the local charts $U_1$ and $U_2$ writes
\begin{equation}\label{main_eq_c2_1_U1}
\begin{array}{rl}
\dot{u}=&c_1+d_1+2c_1u^2-2(b_2-d_2)(u+u^3)+a_1v-u (3a_1u+a_2(-3+u^2))v\\
&-(1+u^2)v^2,\\
\dot{v}=&v\Big(b_2-3b_2u^2+d_2(1+u^2)+a_2v-u(-c_1+d_1+v(2a_1+a_2u+v))\Big),
\end{array}
\end{equation}
and 
\begin{equation}\label{main_eq_c2_1_U2}
\begin{array}{rl}
\dot{u}=&2 b_2(u+u^3)-u(2c_1u+(c_1+d_1)u^3+2d_2(1+u^2))+a_2v-u (3a_2u+a_1(-3+u^2))v\\
&+(1+u^2)v^2,\\
\dot{v}=&-v\Big(b_2+d_2+(c_1+d_1)u^3-a_1v+u^2(-3b_2+d_2+a_1v)+u (c_1+d_1+2a_2v-v^2)\Big),
\end{array}
\end{equation}
respectively.

In the chart $U_2$ system \eqref{main_eq_c2_1_U2} has the origin as an equilibrium point. The Jacobian matrix of system \eqref{main_eq_c2_1_U2} at $(0,0)$ is
$$\left(\begin{matrix}
2b_2-2d_2 & a_2\\
0 & -b_2-d_2
\end{matrix}\right).$$
In order to have a global center of system \eqref{main_eq_c2_1} at the origin, all entries in this Jacobian matrix must be zero. Hence $a_2=b_2=d_2=0.$ Note that then $F=0.$

\smallskip

\noindent{\bf Case} (ii).4:\label{case2.11} 
$c_1\neq 0$. With these new conditions we have that in chart $U_1$ system \eqref{main_eq_c2_1_U1} has two equilibrium points $p_\pm=\left(\pm\sqrt{(-d_1-c_1)/(2c_1)},0\right).$ The Jacobian matrix of system \eqref{main_eq_c2_1_U1} at $p_\pm$ is
$$
\left(\begin{matrix}
\pm2\sqrt{2}\sqrt{-c_1(c_1+d_1)}& \dfrac{a_1(5c_1+3d_1)}{2c_1}\\
0 & \pm(c_1-d_1)\sqrt{\dfrac{-(c_1+d_1)}{2c_1}}
\end{matrix}\right).
$$
Since the Jacobian matrix at any infinite equilibrium point of system \eqref{main_eq_c22_1} must be identically zero, then $a_1=0$ and $d_1=-c_1.$ This implies that the origin is the unique infinite equilibrium point in the chart $U_1.$ To determine the local phase portrait at the origin we must do blow up's. The straight line $u=0$ is not a characteristic direction at the origin of $U_1$ because $R_2(u,v)=-v^3.$ We do the vertical blow-up $(u,v)=(u_1,u_1v_1)$, and system \eqref{main_eq_c2_1_U1} becomes 
\begin{equation*}\label{main_eq_c2_1_U1_blow1_3}
\dot{u_1}=-u_1^2\Big(2d_1+(1+u_1^2)v_1^2\Big),\quad \dot{v_1}=u_1v_1^3.
\end{equation*}
Doing a rescaling of the time we eliminate the common factor $u_1$ between $\dot{u_1}$ and $\dot{v_1}$ and we obtain the system
 \begin{equation}\label{main_eq_c2_1_U1_blow1_resc_3}
\dot{u_1}=-u_1\Big(2d_1+(1+u_1^2)v_1^2\Big),\quad \dot{v_1}=v_1^3.
\end{equation}
The origin is the unique equilibrium point of this system on the straight line $u_1=0.$ The eigenvalues of the Jacobian matrix of system \eqref{main_eq_c2_1_U1_blow1_resc_3} at the origin are $-2d_1$ and $0$. From \cite[Thorem 2.19]{MR2256001} it follows that the origin is a semi-hyperbolic saddle (resp. semi-hyperbolic node) if $d_1>0$ (resp. $d_1<0$, so this case cannot provide a global center). The phase portrait for $d_1>0$ is the one of the left picture of Figure \ref{fig_4}. Going back through the changes of variables we get that the origin of the local chart $U_1$ is formed by two hyperbolic sectors when $d_1>0$. Thus the unique possibility of having a global center in this case is when $d_1>0.$ 

Now we analyze the linearly zero equilibrium at the origin of the chart $U_2$ doing blow up's. The straight line $u=0$ is not a characteristic direction at the origin of $U_2$ because $R_2(u,v)=v(4d_1u^2+v^2).$ Doing the vertical blow up $(u,v)=(u_1,u_1v_1)$ system \eqref{main_eq_c2_1_U2} becomes 
 \begin{equation*}\label{main_eq_c2_1_U2_blow1_31}
\dot{u_1}=u_1^2\Big(2d_1+(1+u_1^2)v_1^2\Big),\quad \dot{v_1}=-u_1\Big(2d_1v_1+v_1^3\Big).
\end{equation*}
Doing a rescaling of the time we eliminate the common factor $u_1$ between $\dot{u_1}$ and $\dot{v_1}$ and we obtain the system
\begin{equation}\label{main_eq_c2_1_U2_blow1_resc_31}
\dot{u_1}=u_1\Big(2d_1+(1+u_1^2)v_1^2\Big),\quad \dot{v_1}=-2d_1v_1-v_1^3.
\end{equation}
The equilibrium points of system \eqref{main_eq_c2_1_U2_blow1_resc_31} on the straight line $u_1=0$ are $(0,0)$ and $q_\pm=\left(0,\pm\sqrt{-2d_1}\right).$ The eigenvalues of the Jacobian matrix of system \eqref{main_eq_c2_1_U2_blow1_resc_31} at the origin are $2d_1$ and $-2d_1$. Since $d_1>0$ the origin is the unique infinite equilibrium point in $U_1$, that is a hyperbolic saddle. Going back through the changes of variables we get that the origin of the local chart $U_2$ is formed by two hyperbolic sectors. The pictures of the blow down in this case coincide with the ones of Figure \ref{fig_4} reversing the sense of all its orbits. 

In summary, for the values of the parameters $a_1=a_2=b_1=b_2=c_2=d_2=0$ and $d_1=-c_1>0.$ the infinite equilibrium points satisfy the conditions of Proposition \ref{prop_main} in order to have a global center. Now we must see if the origin of coordinates is the unique finite equilibrium point. The differential system \eqref{main_eq} reduces to
\begin{equation}\label{bb7}
\dot x=y - 2 c_1 x^2 y, \quad \dot y= -x:
\end{equation}
The unique equilibrium of this system is the $(0,0)$. Therefore by Proposition \ref{prop_main} system \eqref{bb7} has a global center. This completes the proof of statement (g) of Theorem \ref{main_teo}.

\smallskip

\noindent{\bf Case} (ii).5:\label{case2.1} 
$c_1=0$. We assume that $d_1\neq0,$ otherwise system \eqref{main_eq_c2_1} becomes the polynomial differential system of degree two
\begin{equation*}\label{main_eq_c33_1_pc}
\dot{x}=y+2a_1xy,\quad \dot{y}=-x+a_1(x^2-y^2),
\end{equation*}
and we know that such systems has no global centers.
    
These new conditions implies that the origin of the chart $U_2$ is the unique infinite equilibrium point of system \eqref{main_eq_c2_1_U2}, because in the chart $U_1$ system \eqref{main_eq_c2_1_U1} has no infinite equilibrium points. The origin of $U_2$ is a linearly zero equilibrium so we must do blow up's. Since $u=0$ is not a characteristic direction at the origin of $U_2$ because $R_2(u,v)=v(2a_1uv+d_1u^2+v^2)$, we do the vertical blow up $(u,v)=(u_1,u_1v_1).$ Then system \eqref{main_eq_c2_1_U2} becomes 
\begin{equation*}\label{main_eq_c2_1_U2_blow1_0}
\dot{u_1}=-u_1^2\Big(d_1u_1^2-3a_1v_1+a_1u_1^2v_1-v_1^2-u_1^2v_1^2\Big),\quad
\dot{v_1}=-u_1v_1\Big(d_1+v_1(2a_1v_1+v_1)\Big).
\end{equation*}
Doing a rescaling of the time we eliminate the common factor $u_1$ between $\dot{u_1}$ and $\dot{v_1}$ and we obtain the system
\begin{equation}\label{main_eq_c2_1_U2_blow1_resc_0}
\dot{u_1}=-u_1\Big(d_1u_1^2-3a_1v_1+a_1u_1^2v_1-v_1^2-u_1^2v_1^2\Big),\quad
\dot{v_1}=-v_1\Big(d_1+v_1(2a_1v_1+v_1)\Big).
\end{equation}
The equilibrium points of system \eqref{main_eq_c2_1_U2_blow1_resc_0} on the straight line $u_1=0$ are $(0,0)$ and $q_\pm=\left(0,-a_1\pm\sqrt{a_1^2-d_1}\right).$ The eigenvalues of the Jacobian matrix of system \eqref{main_eq_c2_1_U2_blow1_resc_0} at the origin are $0$ and $-d_1$. From \cite[Thorem 2.19]{MR2256001} it follows that the origin is a semi-hyperbolic node and going back through the change of variables we obtain that there are orbits of system \eqref{main_eq_c2_1} going or coming from the infinity. Consequently the center of this system cannot be global. 

\smallskip

\noindent{\bf Case $(iii)$ of Theorem \ref{lemma_main}}.\label{case3main}
From Theorem \ref{lemma_main} we know that $c_1=c_2=0,$ $d_1=-b_1$ and $d_2=b_2.$ Hence system \eqref{main_eq} reduces to the system
\begin{equation}\label{main_eq_c3_1}
\begin{array}{rl}
\dot{x}=&y-2b_2x^3+2a_1xy-4b_1x^2y+2b_2xy^2-a_2(x^2-y^2),\\
\dot{y}=&-x+2a_2xy-2b_2x^2y-4b_1xy^2+2b_2y^3+a_1(x^2-y^2).
\end{array} 
\end{equation}
From \eqref{poincare_comp} system \eqref{main_eq_c3_1} in the local charts $U_1$ and $U_2$ becomes
\begin{equation}\label{main_eq_c3_1_U1}
\begin{array}{rl}
\dot{u}=&-v\Big(a_1(-1+3u^2)+v+u(a_2(-3+u^2)+uv)\Big),\\
\dot{v}=&-v\Big(-4b_1u+2b_2(-1+u^2)-a_2v+uv(2a_1+a_2u+ v)\Big),
\end{array}
\end{equation}
and 
\begin{equation}\label{main_eq_c3_1_U2}
\begin{array}{rl}
\dot{u}=&v\Big(a_2+3a_1u-3a_2u^2+v+u^2(-a_1u+v)\Big),\\
\dot{v}=&v\Big(4b_1u+2b_2(-1+u^2)+a_1v+uv(-2a_2-a_1u+v)\Big),
\end{array}
\end{equation}
respectively.

From systems  \eqref{main_eq_c3_1_U1} and \eqref{main_eq_c3_1_U2} the circle of the infinity is filled with equilibria. Doing a rescaling of the time we eliminate the common factor $v$ between $\dot{u}$ and $\dot{v}$ in the chart $U_1$ and we obtain the system
\begin{equation}\label{main_eq_c3_1_U1_resc}
\begin{array}{rl}
\dot{u}=&-\Big(a_1(-1+3u^2)+v+u(a_2(-3+u^2)+uv)\Big),\\
\dot{v}=&-\Big(-4b_1u+2b_2(-1+u^2)-a_2v+uv(2a_1+a_2u+ v)\Big).
\end{array}
\end{equation}
Now we consider the three subcases.

\noindent\textit{Subcase $1$}: $b_2\neq 0.$ Then $\dot{v}\Big|_{v=0}=2b_2+4b_1u-2b_2u^2$. This means that depending on the sign of $b_2$ there are orbits that go or come from infinity in system \eqref{main_eq_c3_1_U1_resc}, and thus the center of this system cannot global.

\noindent\textit{Subcase $2$}: $b_2=0$ and $b_1\neq0.$ Then $\dot{v}\Big|_{v=0}=4b_1u.$ Using the same argument than before, depending on the sign of $b_1$, there are orbits which go or come from infinity in system \eqref{main_eq_c3_1_U1_resc}, and the center cannot global. 

\noindent\textit{Subcase $3$}: $b_2=0$ and $b_1=0.$ Then system \eqref{main_eq_c3_1} becomes the polynomial system of degree two
\begin{equation*}\label{main_eq_c3_1_pc}
\dot{x}=y+2a_1xy-a_2(x^2-y^2),\quad \dot{y}=-x+2a_2xy+a_1(x^2-y^2),
\end{equation*}
and consequently we cannot have a global center. 

\smallskip

\noindent{\bf Case $(iv)$ of Theorem \ref{lemma_main}}
We have that $c_2=d_1=d_2=0$ and $G=-a_2^2 b_2^3 + 3 a_2^2 b_2b_1^2 + 6 a_2a_1b_2^2 b_1 - 2 a_2a_1b_1^3 + a_1^2 b_2^3 - 3 a_1^2 b_2b_1^2=0.$ Using that $c_2=d_1=d_2=0$ system \eqref{main_eq} becomes
\begin{equation}\label{main_eq_c4_1}
\begin{array}{rl}
\dot{x}=&-b_2 x^3 + y + 2 a_1 x y - (3 b_1 + c_1) x^2 y + 3 b_2 x y^2 -a_2 (x^2 - y^2),\\
\dot{y}=&-x + (b_1 + c_1) x^3 + 2 a_2 x y - 3 b_2 x^2 y + (-3 b_1 + c_1) x y^2 + b_2 y^3+ a_1 (x^2 - y^2).
\end{array} 
\end{equation}
From \eqref{poincare_comp} system \eqref{main_eq_c4_1} in the local charts $U_1$ and $U_2$ writes
\begin{equation*}\label{main_eq_c4_1_U1}
\begin{array}{rl}
\dot{u}=&b_1 + c_1 + 2 c_1 u^2 - 2 b_2 (u + u^3) + a_1 v - v \Big(a_2 u (-3 + u^2)+ v + u^2 (3 a_1 + v)\Big),\\
\dot{v}=&v \Big(b_2 - 3 b_2 u^2 + a_2 v + u \big((3 b_1 + c_1 - v (2 a_1 + a_2 u + v)\big)\Big),
\end{array}
\end{equation*}
and 
\begin{equation}\label{main_eq_c4_1_U2}
\begin{array}{rl}
\dot{u}=&-2 c_1 u^2 - (b_1 + c_1) u^4 + 2 b_2 (u + u^3) + a_2 v -\Big(3 a_2 u+ a_1 (-3 + u^2)\Big)uv \\
&+ (1 + u^2) v^2,\\
\dot{v}=&-v \Big(b_2 - 3 b_2 u^2 - a_1 v + u \big(c_1 + c_1 u^2 + b_1 (-3 + u^2) + 2 a_2 v+ a_1 u v - v^2\big)\Big),
\end{array}
\end{equation}
respectively.

In the chart $U_2$ system \eqref{main_eq_c4_1_U2} has the origin as an equilibrium point. The Jacobian matrix of system \eqref{main_eq_c4_1_U2} at $(0,0)$ is 
$$\left(\begin{matrix}
2b_2 & a_2\\
0 & -b_2
\end{matrix}\right).$$
In order to have a global center of system \eqref{main_eq_c4_1} at the origin all entries in this Jacobian matrix must be zero. So $a_2=b_2=0.$ Note that these new conditions implies that $G=0.$ We divide the study into the following three subcases.

\smallskip

\noindent\textit{Subcase $1$}: $b_1=0$. 
This subcase was already analyzed in case $(i)$ when $b_1= 0$.

\noindent\textit{Subcase $2$}: $b_1\neq0$ and $c_1\neq0$.  
This subcase was already worked in case $(ii).1.$

\noindent\textit{Subcase $3$}: $b_1\neq0$ and $c_1=0$. 
This subcase was already done in case $(ii).2.$

This complete the proof of Theorem \ref{main_teo}.

\section{Acknowledgements}

Jaume Llibre is partially supported by Agencia Estadal de Investigación grant PID2019-104658GB-I00, the H2020 European Research Council grant MSCA-RISE-2017-777911, AGAUR (Generalitat de Catalunya) grant 2021SGR00113, and by the Acadèmia de Ciències i Arts de Barcelona. Gabriel Rondón is supported by São Paulo Research Foundation (FAPESP) grants 2020/06708-9 and 2022/12123-9. 

\bibliographystyle{abbrv}
\bibliography{references1}

\end{document}